\documentclass[a4paper]{article}
\usepackage{graphicx} 
\usepackage{amsmath, amsthm, amsfonts, amssymb, amscd, bm, enumerate}
\usepackage{verbatim}

\usepackage{url}

\newtheorem{theorem}{Theorem}
\newtheorem{corollary}{Corollary}
\newtheorem{lemma}{Lemma}
\newtheorem{definition}{Definition}
\newtheorem{assumption}{Assumption}

\theoremstyle{remark}
\newtheorem{remark}{Remark}
\def\F{\mathcal{F}}
\def\E{\mathcal{E}}

\def\fE{\mathfrak{E}}
\def\D{\mathfrak{D}}

\begin{document}
\title{Representing filtration consistent nonlinear expectations as
$g$-expectations in general probability spaces}
\author{Samuel N. Cohen\\ University of Oxford}

\date{\today}

\maketitle
\begin{abstract}
We consider filtration consistent nonlinear expectations in probability spaces satisfying only the usual conditions and separability. Under a domination assumption, we demonstrate that these nonlinear expectations can be expressed as the solutions to Backward Stochastic Differential Equations with Lipschitz continuous drivers, where both the martingale and the driver terms are permitted to jump, and the martingale representation is infinite dimensional. To establish this result, we show that this domination condition is sufficient to guarantee that the comparison theorem for BSDEs will hold, and we generalise the nonlinear Doob-Meyer decomposition of Peng to a general context.

MSC: 60G48, 60H20, 91B06
\end{abstract}

\section{Introduction}
Much work has been done regarding risk-averse decision making in various contexts. One approach to this has been to assume that agents make decisions based on the `expectation' of a random outcome, but to allow this expectation to be nonlinear. This allows resolution of the famous Allais and Ellsberg paradoxes, while still retaining much of the flavour of classical approaches.

A significant problem in this context is to guarantee that these nonlinear expectations are \emph{time consistent}, that is, that they can be consistently updated using new observations. As many of these nonlinear expectations are not time-consistent, it is useful to give representations for those which are. In \cite{Peng1997} (see also \cite{Peng2004}), Peng gives an axiomatic approach to these nonlinear expectations. In \cite{Peng1997}, of particular interest are the `$g$-expectations', which arise  from the solutions to  Backward Stochastic Differential Equations (BSDEs).

In  Coquet, Hu, M\'emin and Peng \cite{Coquet2002}, it is shown that every nonlinear expectation satisfying a certain domination property must solve a BSDE. At the end of that paper \cite[Remark 7.1]{Coquet2002}, the following comment is made.
\begin{quote}
``In this paper we have limited ourselves to treat the situation where
the filtration is generated by a Brownian motion. A natural question is whether our
nonlinear supermartingale decomposition approach can be applied to more general
situations. A general positive answer seems unlikely, due to the lack of comparison
theorem for BSDE’s driven by discontinuous processes.''
\end{quote}

In this paper, we answer this question in the affirmative, using the BSDEs and comparison theorem in \cite{Cohen2010}. We show that all  nonlinear expectations satisfying a domination property similar to that in \cite{Coquet2002} can be represented by solutions to BSDEs. The domination property which we use is sufficient to guarantee that a comparison theorem holds, and so this extension of \cite{Coquet2002} is possible. We do this making no substantive assumptions on the probability space (we only assume the usual conditions and that $L^2(\F_T)$ is separable). Furthermore, even in the context of a Brownian filtration, our results extend \cite{Coquet2002} to allow a countable number of independent Brownian motions. A weaker extention (to the case of a L\'evy filtration) was obtained by Royer \cite{Royer2006}, and in discrete time there are similar results (see \cite{Cohen2008c}, \cite{Cohen2009a}).

A more general result, restricted to the context of a Brownian filtration, is given by Hu, Ma, Peng and Song \cite{Hu2008}. This result uses a weaker domination property, which corresponds to considering solutions to quadratic BSDEs. As no existence results for quadratic BSDEs are available in the general context considered in \cite{Cohen2010}, we are not yet able to encompass these cases.

Alternative representations exist for nonlinear expectations, for example, Bion-Nadal (\cite{Bion-Nadal2008},\cite{Bion-Nadal2008a}) has a representation for the penalty term of time-consistent convex risk measures (which, up to a change of sign, can be seen to be equivalent to the nonlinear expectations considered here). Similarly, in the Brownian filtration, Delbaen, Peng and Rosazza-Gianin \cite{Delbaen2008} represent these penalty terms using $g$-expectations. The approach of this paper is instead to give a representation of the nonlinear expectation directly, which allows us to avoid any assumption of convexity.

In this paper, we begin by summarizing and generalising the results and approach of \cite{Cohen2010} to BSDEs in general probability spaces. We then also reproduce the key results on filtration-consistent expectations (without proof where the result is exactly as in \cite{Coquet2002}). We proceed to generalise a result of \cite{Peng1999}, giving a Doob-Meyer type decomposition for $g$-expectations in general probabilty spaces, and furthermore, for general nonlinear expectations satisfying our domination property. Finally, using the previous results, we show that any nonlinear expectation satisfying our domination property must equal a $g$-expectation.

\section{BSDEs in General Spaces}

We here give the key results regarding BSDEs in general probability spaces.
These are taken without proof from \cite{Cohen2010}. For simplicity, we shall
restrict our attention to the scalar case. As usual, unless otherwise indicated, all (in-)equalities should be read as `up
to evanescence'.

\begin{assumption}\label{assn1}
 We shall henceforth assume that
\begin{enumerate}[(i)]
 \item the usual conditions hold on our filtered probability space $(\Omega, \F, \{\F_t\}_{t\in[0,T]}, \mathbb{P})$, and $\mathcal{F}_0$ is the $\mathbb{P}$-completion of the trivial $\sigma$-algebra $\{\Omega, \emptyset\}$,
 \item $L^2(\mathcal{F}_T)$ is separable, and
 \item we have some (arbitrary,) deterministic, strictly increasing process $\mu$ with $\mu_T<\infty$.
\end{enumerate}
\end{assumption}

\begin{remark}
The process $\mu$ will be used in the place of Lebesgue measure in our
BSDE. The assumption that $\F_0$ is trivial is not strictly necessary, but is used to simplify notation (as it implies no martingale has a jump at zero).
\end{remark}

\begin{definition} \label{defn:inducedmeasure}
For any nondecreasing process of finite variation $\mu$, we define the measure
induced by $\mu$ to be the measure over $\Omega\times [0,T]$ given by
\[A\mapsto E\left[\int_{[0,T]}I_A(\omega, t) d\mu\right].\]
Here $A\in\mathcal{P}$, the predictable $\sigma$-algebra, and the integral is
taken pathwise in a Stieltjes sense.
\end{definition}

The following version of the martingale representation theorem (from Davis and
Varaiya \cite{Davis1974}, see also Kunita and Watanabe \cite{Kunita1967} and Malamud \cite{Malamud2007}) is fundamental to our approach.

\begin{theorem}[Martingale Representation Theorem; \cite{Davis1974}]
\label{thm:martrep}
Suppose $L^2(\mathcal{F}_T)$ is a separable Hilbert space, with an inner product
$(X,Y) = E[XY]$. Then there exists a finite or countable sequence of
square-integrable $\{\mathcal{F}_t\}$-martingales $M^1, M^2,...$ such that every
square integrable $\{\mathcal{F}_t\}$-martingale $N$ has a representation
\[N_t= N_0+\sum_{i=1}^{\infty} \int_{]0,t]}Z^i_u dM^i_u\]
for some sequence of predictable processes $Z^i$.  This sequence satisfies
\begin{equation} \label{eq:Zbound}
E\left[\sum_{i=0}^\infty\int_{]0,T]} (Z^i_u)^2 d\langle M^i\rangle_u\right]
<+\infty.
\end{equation}

These martingales are orthogonal (that is, $E[M^i_TM^j_T]=0$ for all $i\neq j$),
and the predictable quadratic variation processes $\langle M^i\rangle$ satisfy
\[\langle M^1\rangle \succ \langle M^2\rangle \succ \ldots,\]
where $\succ$ denotes absolute continuity of the induced measures (Definition
\ref{defn:inducedmeasure}). Furthermore, these martingales are unique, in that
if $N^i$ is another such sequence, then $\langle N^i\rangle \sim \langle
M^i\rangle$, where $\sim$ denotes equivalence of the induced measures.
\end{theorem}

We shall denote by $\mathbb{R}^\infty$ the set of countable sequences of real
values.

\begin{definition}\label{defn:Mnorm}
We define the stochastic seminorm $\|\cdot\|_{M_t}$ on $\mathbb{R}^{\infty}$ as follows. For each $i\in \mathbb{N}$, consider $\langle M^i\rangle$ as a measure on the predictable $\sigma$-algebra. Let $\langle M^i\rangle$ have the Lebesgue-decomposition
\[\langle M^i\rangle_t = m^{i,1}_t + m^{i,2}_t,\]
where $m^{i,1}_t$ is absolutely continuous with respect to $\mu\times \mathbb{P}$ and $m^{i,2}_t$ is orthogonal to $\mu\times \mathbb{P}$. As they represent bounded measures on the predictable $\sigma$-algebra, both $m^{i,1}_t$ and $m^{i,2}_t$ will be nondecreasing predictable processes. As measures, we can find a version $\phi^i$ of the Radon-Nikodym derivative
\[\phi^i=\frac{dm^{i,1}}{d(\mu\times \mathbb{P})}\]
such that $\phi^i = 0$, $m^{i,2}$-a.e.

We define, for $z_t\in \mathbb{R}^{\infty}$,
\[\|z_t\|^2_{M_t} := \sum_i \left[|z_t^i|^2 \phi^i(t,\omega)\right]\]
where $z^i_t\in\mathbb{R}$ is the $i$'th element in $z_t$. 

We note that, for any predictable, progressively measurable process $Z$ taking values in $\mathbb{R}^{\infty}$, and in particular for processes satisfying (\ref{eq:Zbound}), we have the inequality
\begin{equation}\label{eq:isometry}
\begin{split}
&E\left[\int_{A}\|Z_t\|^2_{M_t} d\mu\right] \leq E\left[\sum_i\int_{A}(Z_t^i)^2 d\langle M^i_t\rangle\right]\\
&=E\left[\sum_i\left(\int_{A}Z_t^idM^i_t\right)^2\right]=E\left[\left(\sum_i\int_{A}Z_t^idM^i_t\right)^2\right] \end{split}
\end{equation}
for any predictable set $A\subseteq \Omega\times[0,T]$. (Note the latter equalities are simply the standard isometry used in the construction of the stochastic integral, by the orthogonality of the $M^i$.)

For any predictable process $Z$ taking values in $\mathbb{R}^\infty$ with (\ref{eq:isometry}) finite, any
predictable set $A$, for notational simplicity we shall write
\[\begin{split}
 \int_A Z_t dM_t &:= \sum_i \int_A Z^i_t dM^i_t,\\
 Z_t\Delta M_t &:= \sum_i Z_t^i \Delta M^i_t,\\
 \int_A Z_t^2 d\langle M\rangle_t &:= \sum_i \int_A |Z^i_t|^2 d\langle
M^i\rangle_t
\end{split}\]
\end{definition}

\begin{definition} \label{defn:processspaces}
We define the following spaces
\[\begin{split}
H^2_M &=\left\{Z:\Omega\times[0,T]\to\mathbb{R}^{\infty},
\quad\text{predictable, } E\left[\int_{]0,T]}Z_t^2 d\langle
M\rangle_t\right] <+\infty\right\},\\
S^2 &= \left\{Y:\Omega\times[0,T]\to\mathbb{R}, \quad\text{adapted, }
E\left[\sup_{t\in[0,T]} \|Y_t\|^2\right]<+\infty\right\},\\
H^2_\mu&= \left\{Y:\Omega\times [0,T] \to\mathbb{R}, \quad\text{progressive,
}\int_{]0,T]} E[\|Y_{t}\|^2] d\mu_t<+\infty\right\},
\end{split}\]
where two elements $Z, \bar Z$ of $H^2_M$ are deemed equivalent if
\[E\left[\int_{[0,T]}(Z_t-\bar Z_t)^2d\langle M\rangle_t\right]
=0,\]
two elements of $S^2$ are deemed equivalent if they are indistinguishable, and
two elements of $H^2_\mu$ are equivalent if they are equal
$\mu\times\mathbb{P}$-a.s.
\end{definition}

\begin{remark}
We note that $H^2_M$ is itself a complete metric space, with norm given by 
$Z\mapsto E\left[\int_{]0,T]}Z_t^2d\langle M\rangle_t\right]$.
Similarly for $H^2_\mu$. Note also that the martingale representations
constructed in Theorem \ref{thm:martrep} are unique in $H^2_M$.
\end{remark}

\begin{theorem} \label{thm:BSDEExist2}
 Let $g:\Omega\times[0,T]\times
\mathbb{R}\times\mathbb{R}^{\infty}\to\mathbb{R}$ be a predictable function such
that
\begin{itemize}
\item $E\left[\int_{]0,T]} |g(\omega, t, 0,\mathbf{0})|^2d\mu_t\right]<+\infty$
\item There exists a quadratic firm Lipschitz bound on $F$, that is, a
measurable deterministic function $c_t$ uniformly bounded by some
$c\in\mathbb{R}$, such that, for all $y, y'\in \mathbb{R}$, $z,
z'\in\mathbb{R}^{\infty}$, all $t>0$
\[\|g(\omega, t, y, z) - g(\omega, t, y', z')\|^2\leq c_t
\|y_t-y'\|^2+ c\|z-z'\|^2_{M_t} \qquad d\mu\times d\mathbb{P}-a.s.\]
and
\[c_t(\Delta \mu_t)^2<1 \qquad \text{ for all }t>0.\]
 Note that the variable bound $c_t$ need only apply to the behaviour of $F$ with
respect to $y$.
\end{itemize}
A function satisfying these conditions will be called \emph{standard}.
Then for any $Q\in L^2(\mathcal{F}_T)$, the BSDE
with driver $g$
\begin{equation} \label{eq:BSDE}
 Y_t - \int_{]t,T]} g(\omega, u, Y_{u-}, Z_u) d\mu_u + \int_{]t,T]} Z_u
dM_u = Q
\end{equation}
has a unique solution $(Y,Z)\in S^2\times H^2_M$.
\end{theorem}

From this point onwards, for notational simplicity, we shall regard $\omega$ as implicit in the function $g$, whenever this does not lead to confusion.
\begin{remark}\label{rem:propertiesofg}
 Note that the behaviour of $g$ at $t=0$ is irrelevant to the solution of the BSDE, however we still obtain a solution with values $(Y_t, Z_t)$ for $t\in[0,T]$. Note also that for any $y\in \mathbb{R}, z,z'\in\mathbb{R}^\infty$, we know $\|g(t, y, z) - g(t, y, z')\|=0$ $m^{i,2}$-a.e. for all $i$, by the definition of the norm $\|\cdot\|_{M_t}$.
\end{remark}

\begin{theorem}[Comparison Theorem]\label{thm:CompThm}
Suppose we have two BSDEs corresponding to standard coefficients and terminal
values $(g, Q)$ and $(g', Q')$ . Let $(Y, Z)$ and $(Y', Z')$ be
the associated solutions. Suppose that for some $s$, the following conditions
hold:
\begin{enumerate}[(i)]
	\item $Q\geq Q'$ $\mathbb{P}$-a.s.
	\item $\mu \times \mathbb{P}$-a.s. on $[s,T]\times \Omega$,
\[g(u, Y'_{u-}, Z'_u) \geq g'(u, Y'_{u-}, Z'_u).\]
	\item There exists a measure $\tilde{\mathbb{P}}$
equivalent to $\mathbb{P}$ such that
	\[X_r := -\int_{]s,r]} [g(\omega, u, Y'_{u-}, Z_u) - g(u, Y'_{u-}, Z'_u)]d\mu_u +
\int_{]s,r]}[Z_u -Z'_u]dM_u\]
	is a $\tilde{\mathbb{P}}$ supermartingale on $[s,T]$.
\end{enumerate}
 It is then true that $Y \geq Y'$ on $[s,T]\times \Omega$, except possibly on some
evanescent set. Furthermore, this comparison is strict, that is, for any $s$ and
any $A\in\mathcal{F}_s$ such that $Y_s=Y'_s$ $\mathbb{P}$-a.s. on $A$, we
have $Y_u = Y'_u$ on $[s,T]\times A$, up to evanescence.
\end{theorem}

\begin{definition}
In light of this, we make the following definition.

If $g$ is such that condition (iii) of Theorem \ref{thm:CompThm} holds for any
special semimartingales $Y, Y'\in H^2_\mu$, (where $Z$ and $Z'$ are from
the martingale representation theorem applied to the martingale parts of $Y$ and
$Y'$) then $g$ shall be called \emph{balanced}.
\end{definition}

\begin{lemma}\label{lem:wellbalancedisbalanced}
If for any $y\in \mathbb{R}$, any $z,z' \in \mathbb{R}^{\infty}$
\[\frac{|g(t, y, z)-g(t,y,z')|}{\|z-\bar z\|^2_{M_t}}
|(z-z') \Delta M_t| <1\]
up to evanescence, then $g$ is balanced.
\end{lemma}

To prove this, we first need the following lemma, based on results of Lepingle and M\'emin \cite{Lepingle1978}, (see also Protter and Shimbo, \cite{Protter2008}).

\begin{definition}[Dol\'eans-Dade Exponential]
 Let $N$ be a local martingale. Then we shall write
\[\mathfrak{E}(N;t) := \exp(N_t - \langle N^c\rangle_t/2) \prod_{0<s\leq t} (1+\Delta N_s)\exp(-\Delta N_s),\]
which is the solution $\mathfrak{E}(N;t) = M_t$ of the equation
\[M_t = 1+\int_{]0,t]} M_{s-} dN_s.\]
\end{definition}

\begin{lemma}\label{lem:exponentialmoments}
 Let $N$ be a square-integrable martingale, with $\langle N\rangle$ bounded. Then $\fE(N;\cdot)$ is a martingale, and for any $p>0$, $E[|\mathfrak{E}(N;T)|^p]<\infty$.
\end{lemma}
\begin{proof}
It is clear that $\fE(N;\cdot)$ is a local martingale, by Lepingle and M\'emin \cite[Thm II.2]{Lepingle1978} it is a square integrable martingale. It is easy to verify that
\[\mathfrak{E}^2(N;t) = 1+ \int_{]0,t]} \mathfrak{E}^2(N;s-)d(2N+[N])_s=\mathfrak{E}(2N+[N];t).\]
As $\langle N\rangle\leq k$ for some $k$, we can write
\begin{equation}\label{eq:expsquarebound}
\mathfrak{E}^2(N;t) = \mathfrak{E}(2N+[N]-\langle N\rangle +\langle N\rangle;t)\leq e^{k}\mathfrak{E}(2N+[N]-\langle N\rangle;t). 
\end{equation}
We now see that $\tilde N := 2N+[N]-\langle N\rangle = 2N +[N^d]-\langle N^d\rangle$ and this is a local martingale, hence
\[\langle \tilde N^c\rangle = 2\langle N^c\rangle\leq 2k\]
and
\[(\Delta\tilde N)^2 = (3\Delta N-\Delta\langle N^d\rangle))^2\leq 18(\Delta N)^2+2(\Delta\langle N^d\rangle)^2.\]
These quantities are integrable at $T$, so $\tilde N$ is a square-integrable martingale. Furthermore,
\[\langle \tilde N^d \rangle \leq 18 \langle N^d\rangle + 2\sum_{0< u\leq t}((\Delta\langle N^d\rangle)^2) \leq 18 \langle N^d\rangle + 2\langle N^d\rangle^2\leq 18k+2k^2,\]
and we see that $\langle \tilde N \rangle \leq 20k+2k^2$, in particular, that this is a finite bound. Hence $\tilde N$ is a square-integrable martingale with $\langle \tilde N \rangle$ bounded.

From \cite[Thm II.2]{Lepingle1978}, we see that $\mathfrak{E}(\tilde N;t)$ is a square integrable martingale, and from (\ref{eq:expsquarebound})
\[E[(\mathfrak{E}(N;T))^4]\leq e^{2k} E[(\mathfrak{E}(\tilde N;T))^2]<\infty.\]

We now iterate this process, noticing that $\tilde N$ satisfies the requirements of the lemma, and hence if $\tilde{\tilde N} = 2\tilde N +[\tilde N]-\langle \tilde N\rangle$, (which is, by the same logic, a square integrable martingale with $\langle \tilde{\tilde N}\rangle$ bounded), 
\[E[(\mathfrak{E}(N;T))^8] = E[(\mathfrak{E}(\tilde N;T))^4]\leq e^{2(20k+2k^2)} E[(\mathfrak{E}(\tilde{\tilde N};T))^2]<\infty.\]
Hence we obtain, after $n$ iterations,
\[E[(\mathfrak{E}(N;T))^{2^n}]<\infty\]
and by Jensen's inequality, the result is proven for any finite $p$.
\end{proof}

\begin{proof}[Proof of Lemma \ref{lem:wellbalancedisbalanced}]
Define
\[N_t = \int_{]0,t]} \left(\frac{g( u, Y'_u,Z_u)-g(u,Y'_u,Z'_u)}{\|Z_u-Z'_u\|^2_{M_u}}\right)(Z_u-Z'_u) dM_u\]
Let $\Lambda$ be the process defined by the Dol\'eans-Dade exponential
\[\Lambda_t = 1+ \int_{]0,t]} \Lambda_{u-}dN_u = \mathfrak{E}(N;t).\]
By the assumption of the Lemma, we see that $|\Delta N_t| < 1$, and
so $\Lambda_t$ is a strictly positive local Martingale.
Furthermore, we know that $N$ has predictable quadratic variation
\[\begin{split}
\langle N\rangle_t &= \int_{]0,t]} \left(\frac{g(u,Y'_u,Z_u)-g(u,Y'_u,Z'_u)}{\|Z_u-Z'_u\|^2_{M_u}}\right)^2(Z_u-Z'_u)^2 d\langle M\rangle_u \\
&= \int_{]0,t]} \frac{(g( u,
Y'_u,Z_u)-g(u,Y'_u,Z'_u))^2}{\|Z_u-Z'_u\|^2_{M_u}}d\mu_u
\\
&\leq c\mu_t
\end{split}\]
where $c$ is the Lipschitz constant of $g$, using the decomposition $d\langle M^i\rangle = \phi^i d\mu + dm^{i,2}$ and Remark \ref{rem:propertiesofg}. By Lemma \ref{lem:exponentialmoments}, this shows that $\Lambda$ has moments of all orders, and is a true martingale on $[0,T]$. We can therefore define the measure
 $\tilde{\mathbb{P}}$ by $d\tilde{\mathbb{P}}/d\mathbb{P}=\Lambda_T$.

 By Girsanov's theorem (see \cite[Theorem 3.11]{Jacod2003}), we see that
 \[\tilde M^i_t = M_t^i - \int_{]0,t]} \frac{g(\omega, u, Y'_u,Z_u)-g(u,Y'_u,Z'_u)}{\|Z_u-Z'_u\|^2_{M_u}}(Z_u-Z'_u)^id\langle M^i\rangle_u\]
 is a $\tilde{\mathbb{P}}$-local martingale. Hence
 \[\begin{split}
    X_t &= \int_{]0,t]} (Z_u-Z'_u)d\tilde M_u\\
 &=-\int_{]0,t]} (g( u, Y'_u,Z_u)-g(u,Y'_u,Z'_u)) d\mu_u + \int_{]0,t]} (Z_u- Z'_u)dM_u
   \end{split}\]
is a $\tilde{\mathbb{P}}$-local martingale.

Finally, by H\"older's inequality, for any stopping time $\tau$, any $\epsilon\in ]0,2]$
\[E_{\tilde{\mathbb{P}}}[X^{2-\epsilon}_\tau] = E_{\mathbb{P}}[\Lambda_T X_{\tau}^{2-\epsilon}] \leq E_{\mathbb{P}}[\Lambda_T^{2/\epsilon}]^{(\epsilon/2)} E_{\mathbb{P}}[X_{\tau}^{2}]^{1-\epsilon/2},\]
which is uniformly bounded, by Lemma \ref{lem:exponentialmoments} and the fact $X$ is $\mathbb{P}$-square-integrable. It follows that $X$ is a true $\tilde{\mathbb{P}}$-martingale.
\end{proof}

\subsection{A scalar extension}
As we are considering the case of scalar-valued BSDEs, it is useful to extend our existence result beyond the firmly Lipschitz assumptions of \cite{Cohen2010}, as this will enable us to use various penalisation methods.
\begin{theorem}\label{thm:BSDEExist3}
 Let $g:\Omega\times[0,T]\times \mathbb{R}\times \mathbb{R}^\infty$ be a predictable function such that
\begin{enumerate}
 \item $E\left[\int_{]0,T]}g(t, 0,\mathbf{0})^2d\mu_t\right] <+\infty$
\item $g$ is Lipschitz, that is, there exists $c\in\mathbb{R}$ such that for any $y, y' \in \mathbb{R}$, any $z,z'\in \mathbb{R}^\infty$
\[\|g(t, y,z)-g(t, y, z)\|^2\leq c(\|y-y'\|^2 + \|z-z'\|^2_{M_t})\qquad d\mu\times d\mathbb{P}-a.s.\]
and furthermore, for all $y\neq y'$, $g$ satisfies
\[\left(\frac{g(t,y,z)-g(t,y',z)}{y-y'}\right)\Delta\mu_t \leq 1-(1+c)^{-1}.\]
\end{enumerate}
Then for any $Q\in L^2(\mathcal{F}_T)$, the BSDE with driver $g$ has a unique solution $(Y,Z)\in S^2\times H^2_M$. Furthermore, if $g$ is balanced (that is, condition (iii) of Theorem \ref{thm:CompThm} is satisfied), then the comparison theorem holds.
\end{theorem}

\begin{proof}
 As $g$ is Lipschitz with constant $c$ and $\mu$ is a finite valued increasing process, there are at most finitely many times $t_1, t_2,..., t_k$ such that $c(\Delta\mu_t)^2\geq 1$ (and these times are deterministic). Hence, between these times, we have a standard BSDE. We shall show that,
\begin{enumerate}[(i)]
 \item For each $t_i$, we can take any $Y_{t_i} \in L^2(\mathcal{F}_{t_i})$, and obtain a unique pair $(Y_{t_i*}, Z_{t_i})$, where $Y_{t_i*}\in L^2(\mathcal{F}_{t_i-})$, $Z_{t_i}$ is $\mathcal{F}_{t_i-}$-measurable and $Z_{t_i}\Delta M_{t_i} \in L^2(\mathcal{F}_{t_i})$.
\item We can then use this value $Y_{t_i*}$ as the terminal value for a BSDE on the interval $[t_{i-1}, t_{i}[$, which has a unique solution, as our driver is standard (recalling that the behaviour of the driver at the left-endpoint is unimportant for the BSDE solution).
\item The BSDEs we construct on $[t_{i-1}, t_{i}[$ satisfy $\lim_{t\uparrow t_i}Y_t = Y_{t_i*}$ almost surely, so our solutions satisfy $Y_{t_i*} = Y_{t_i-}$ up to evanescence.
\end{enumerate}
 Backward induction then yields that we have a solution to the BSDE on $[0,T]$. Note that, as $\{t_1,...,t_k\}$ is finite, the processes we construct are appropriately predictable.

We first show that (i) our solution can be constructed at each problematic jump-time $t_i$. At $t_i$, we have the equation
\[Y_{t_i} = Y_{t_i*} - g(t_i, Y_{t_i*}, Z_{t_i}) \Delta \mu_{t_i} + Z_{t_i} \Delta M_{t_i},\]
where $(Y_{t_i*}, Z_{t_i})$ are to be determined. Taking an expectation and difference, we see that $Z_{t_i} \Delta M_{t_i} = Y_{t_i} - E[Y_{t_i}|\mathcal{F}_{t_i-}]$. As this is a martingale difference, by the martingale representation theorem, we obtain a solution $Z_{t_i}$. Fixing $Z_{t_i}$ at this solution, we then see that
\[E[Y_{t_i}|\mathcal{F}_{t_i-}] = Y_{t_i*} - g(t_i, Y_{t_i*}, Z_{t_i}) \Delta \mu_{t_i}.\]
Writing $\phi(y) :=y-g(t_i, y, Z_{t_i})\Delta \mu_{t_i}$, our assumptions on $g$ show that $\phi$ is bi-Lipschitz with constant $(1+c)$, and strictly increasing. Hence it has a strictly increasing bi-Lipschitz inverse, also with constant $(1+c)$. We therefore define $Y_{t_i*} = \phi^{-1}(E[Y_{t_i}|\mathcal{F}_{t_i-}])$. By Lipschitz continuity and Jensen's inequality, $Y_{t_i*}\in L^2(\mathcal{F}_{t_i-})$.

We now consider (ii), our BSDE on an interval $]t_{i-1},t_i[$. As $g$ is standard on this interval, $g':=g(t,y,z)I_{t\neq t_i}$ is standard on $]t_{i-1}, t_i]$. Hence it has a solution $(Y', Z')$ on $[t_{i-1},t_i]$, with $Y'_{t_i}=Y_{t_i*}$. As we have a terminal value which is $\mathcal{F}_{t_i-}$-measurable, it is easy to verify that our solution will satisfy $Z'_{t_i}\equiv0$. We see that this is idential to the BSDE with driver $g$ written on the interval $]t_{i-1},t_i[$, and so we can define our solution $(Y_{t}, Z_t) = (Y'_t, Z'_t)$ for all $t\in [t_{i-1},t_i[$. Note that as $Z'_{t_i}\equiv 0$ and $g'(t_i, \cdot, \cdot) \equiv 0$, we also have (iii), $Y'_{t_i-} = Y'_{t_i} = Y_{t_i*}$.

For the comparison theorem, we immediately see that it holds on each interval $[t_{i-1}, t_i[$. At $t_i$, we have an essentially identical argument as that given in discrete time in \cite[Theorems 3.2 and 3.5]{Cohen2009a}.
\end{proof}

\begin{remark}
 Note that, if $g$ is Lipschitz continuous and nonincreasing in $y$, then it is easy to verify that condition (2) holds.
\end{remark}

\subsection{Gr\"onwall's inequality}
In \cite{Cohen2010}, we also derive a version of Gr\"onwall's inequality, which
shall be useful here. 

\begin{definition}\label{defn:jumpinversion}
Let $\nu$ be a c\`adl\`ag function of finite variation with $\Delta\nu_t<1$ for all $t$. The \emph{right-jump-inversion} of $\nu$ is defined by
\[\tilde \nu_t := \nu_t + \sum_{0\leq s\leq t} \frac{(\Delta
\nu_s)^2}{1-\Delta\nu_s}.\]
And satisfies $\fE(-\nu;t)= \fE(\tilde \nu; t)^{-1}$.
\end{definition}

\begin{definition}
Let $u, v$ be two measures on a $\sigma$-algebra $\mathcal{A}$. We write $du\leq
dv$ if, for any $A\in \mathcal{A}$, $u(A) \leq v(A)$.
\end{definition}

\begin{lemma}[Backward Gr\"onwall Inequality] \label{lem:groenwallinequality}
Let $u$ be a process such that, for $\nu$ a nonnegative Stieltjes measure with
$\Delta\nu_t<1$ and $\alpha$ a $\tilde\nu$-integrable process, $u$ is
$\nu$-integrable and
\[u_{t} \leq \alpha_t + \int_{]t,T]} u_{s} d\nu_s,\]
then
\[u_{t} \leq \alpha_t + \mathfrak{E}(-\nu_t)\int_{]t,T]}\mathfrak{E}(\tilde
\nu_{s}) \alpha_s d\tilde\nu_s.\]
If $\alpha_t=\alpha$ is constant, this simplifies to
\[u_{t} \leq \alpha \mathfrak{E}(\tilde \nu;T)
\mathfrak{E}(\tilde\nu;t)^{-1}=\alpha\mathfrak{E}(-\nu;t)\mathfrak{E}(-\nu;T)^{
-1}.\]
\end{lemma}

\section{Filtration Consistent Expectations}
We now reproduce, for completeness, relevant results from Coquet et al
\cite{Coquet2002}. These are given without proof where the argument of
\cite{Coquet2002} carries over without change, or is standard.

\begin{definition}\label{defn:Fexpectations}
 A nonlinear expectation is a functional $\E: L^2(\F_T) \to \mathbb{R}$ which
satisfies \emph{strict monotonicity}:
\[\begin{split}
  &\text{if } Q\geq Q' \text{ then }  \E(Q)\geq\E(Q'), \text{
and}\\
   &\text{if } Q\geq Q' \text{ and } \E(Q)=\E(Q') \text{ then }
Q=Q'.
  \end{split}\]
and preserving of constants: $\E(c)=c$ for all constants $c$.

A nonlinear expectation is \emph{filtration consistent} (or $\F$-consistent) if for
each $Q\in L^2(\F_T)$ and each $t\in [0,T]$ there exists a random variable
$\eta\in L^2(\F_t)$ such that $\E(I_A Q) = \E(I_A \eta)$ for all $A\in\F_t$. Such a nonlinear expectation is called an $\F$-expectation.

The following lemma proves that $\eta$ is unique. It is denoted 
$\E(Q|\F_t)$, and is called the conditional $\F$-expectation of $Q$ under
$\F_t$.

An $\F$-expectation $\E$ will be called \emph{translation invariant} if $\E(Q+q|\F_t)=\E(Q|\F_t)+q$ for all $q\in L^2(\F_t)$, all $Q\in L^2(\F_T)$. It is called \emph{convex} if, for any $Q, Q' \in L^2(\F_T)$, any $\lambda\in[0,1]$, $\E(\lambda Q + (1-\lambda) Q') \leq \lambda \E(Q) + (1-\lambda)\E(Q')$. It is called \emph{positively homogenous} if, for any $Q\in L^2(\F_T)$, any $\lambda\geq 0$, $\E(\lambda Q) = \lambda \E(Q)$.
\end{definition}

\begin{lemma}\label{lem:condexpisunique}
 Let $t\leq T$ and $\eta_1, \eta_2 \in L^2(\F_t)$. If $\E(\eta_1 I_A) =
\E(\eta_2 I_A)$ for all $A\in \F_t$, then $\eta_1=\eta_2$.
\end{lemma}

\begin{lemma}\label{lem:regularityofEt}
 Let $\E$ be an $\F$-expectation. Then the following properties hold for all
$Q,  Q'\in L^2(\F_T)$.
\begin{enumerate}[(i)]
 \item For each $0\leq s\leq t\leq T$, $\E(\E(Q|F_t)|F_s)=\E(Q|F_s)$, and in
particular, $\E(\E(Q|\F_t))=\E(Q)$.
\item For any $t$, for all $A\in\F_t$, $\E(Q I_A|\F_t) = I_A \E(Q|\F_t)$.
\item For any $t$, for all $A\in\F_t$, $\E(QI_A + Q' I_{A^c}|\F_t) = \E(QI_A
|\F_t)+ \E(\bar Q I_{A^c}|\F_t)$.
\item For any $t$, if $Q\geq \bar Q$, then $\E(Q|\F_t) \geq \E(Q'|\F_t)$. If
moreover $\E(Q|\F_t) \geq \E(Q'|\F_t)$ for some $t$, then $Q=Q'$.
\end{enumerate}
\end{lemma}

\begin{definition}
 For a given $\F$-expectation $\E$, a process $Y \in S^2$ is called an $\E$-supermartingale if $Y_s\geq\E(Y_t|\F_s)$ a.s. for all $s\leq t$. Similarly, $Y$ is an $\E$-submartingale if $Y_s\leq\E(Y_t|\F_s)$, and an $\E$-martingale if $Y_s=\E(Y_t|\F_s)$.
\end{definition}

\begin{lemma}\label{lem:negsupissub}
 If $\E$ is convex and $Y$ is an $\E$-supermartingale, then $-Y$ is an $\E$-submartingale. 
\end{lemma}

\begin{lemma}
 If $g$ is a balanced driver, and is convex (resp. positively homogenous), then $\E_g$ is convex (resp. positively honogenous).
\end{lemma}

\begin{lemma}\label{lem:sumsupissup}
 If $\E$ is convex and positively homogenous, then the sum of two $\E$-supermartingales is an $\E$-supermartingale.
\end{lemma}

\begin{theorem}[Up/Downcrossing inequalities]\label{thm:updowncrossing}
Let $\E$ be a convex, translation invariant and positively homogenous $\F$-expectation, and $Y$ be an $\E$-submartingale. For any stopping time $S\leq T$, let $M(\omega, Y^S;[\alpha,\beta])$ (resp. $D(\omega, Y^S;[\alpha,\beta])$) denote the number of upcrossings (resp. downcrossings) of the interval $[\alpha, \beta]$ by $Y$ on the interval $[0,S]$.

Then
\[\begin{split}
   \E(M(\omega, Y^S;[\alpha, \beta])) &\leq (\beta-\alpha)^{-1}
(\E((Y_S-\alpha)^+)-(Y_0-\alpha)^+)\\
\E(D(\omega, Y^S;[\alpha, \beta])) &\leq -(\beta-\alpha)^{-1}
\E(-(Y_S-\beta)^+) \\
&\leq (\beta-\alpha)^{-1} \E((Y_S-\beta)^+)
  \end{split}
\]
\end{theorem}
\begin{proof}
 See \cite{Cohen2011}.
\end{proof}
We shall use this result to prove the existence of c\`adl\`ag modifications to nonlinear martingales, see Theorem \ref{thm:cadlagmodn}.

\subsection{$g$-expectations}
\begin{theorem}\label{thm:BSDEisFexp}
Let $g$ be a balanced driver which satisfies
\begin{equation}\label{eq:gzero}
 g(\omega, t, y, \mathbf{0}) =0,\qquad \mu\times
\mathbb{P}
-a.s.
\end{equation}
 Then the operator defined by
\[\E_g(Q|\mathcal{F}_t) := Y_t\]
where $Y$ is the solution to a BSDE (\ref{eq:BSDE}) with driver $g$, is a
conditional $\F$-expectation. $\E_g$ is called the $g$-expectation.
\end{theorem}

\begin{lemma}
If a balanced driver $g(t, z)$ does not depend on $y$, then the $g$-expectation is translation invariant.
\end{lemma}

\begin{lemma}\label{lem:expectationbound}
 Let $g$ be as in Theorem \ref{thm:BSDEisFexp}, and be balanced. Then for
every real $\epsilon >0$, there exists a constant $C_\epsilon$ such that for
every $Q\in L^{2\vee(1+\epsilon)}(\F_T)$,
\[|\E_g(Q)| \leq C_\epsilon \|Q\|_{1+\epsilon},\]
Where $\|\cdot\|_{1+\epsilon}$ is the standard norm in $L^{1+\epsilon}(\F_T)$.
\end{lemma}
\begin{proof}
Define the measure
\[\frac{d\tilde{\mathbb{P}}}{d\mathbb{P}} = \Lambda_T = \mathfrak{E}\left(\int_{]0,T]} \frac{g(s, Y_{s-}, Z_s)}{\|Z_s\|_{M_s}^2} Z_s dM_s; T\right).\]
Similarly to in the proof of Lemma \ref{lem:wellbalancedisbalanced}, this is a stochastic exponential of the form considered in Lemma \ref{lem:exponentialmoments}. Hence $\tilde{\mathbb{P}}$ is a probability measure and $\Lambda_T$ has finite $p$th moment, for any $p$. By Girsanov's theorem, $\E_g(Q)=E_{\tilde{\mathbb{P}}}[Q] = E[\Lambda_T Q]$. By H\"older's inequality, we have $|\E_g(Q)| \leq \|\Lambda_T\|_{1+\epsilon^{-1}} \|Q\|_{1+\epsilon}$, and the claim follows.
\end{proof}

\subsection{$\E^r$ expectations}
We now consider a particularly useful class of $g$-expectations, which we call $\E^r$-expectations.

\begin{definition}
Let $r$ be a predictable process taking values in the space of real-valued countable dimensional matrices $\mathbb{R}^{\infty\times\infty}$, that is, $(r^{ij}) \in \mathbb{R}$ for all $i,j\in\mathbb{N}$.

We denote by $r_tz_t$ the vector in $\mathbb{R}^\infty$ with values $(r_tz_t)^{i} = \sum_j r^{ij}z^j$. (If $z$ were thought of as a column vector, then this would correspond to the classical matrix-vector product.)

The map $z\mapsto rz$ is a linear operator on $H^2_M$. We suppose that  $r$ is uniformly bounded in a modified operator norm, which we denote $\|\cdot\|_{D_t}$, that is, there is $c\in\mathbb{R}$ such that, up to evanescence
\[\|r_t\|_{D_t}^2: = \sup_{z\in H^2_M}\left\{\frac{\|r_tz_t\|_{M_t}^2}{\|z_t\|_{M_t}^2}\right\} = \sup_{\{u\in \mathbb{R}^\infty: \|u\|_{M_t}=1\}}\{\|r_tu\|_{M_t}^2\} <c.\]

The process $r$ will be called uniformly balanced if
\[ \|r_tu\|_{M_t}\times |u\Delta M| < 1\]
for all $u\in \mathbb{R}^\infty$ with $\|u\|_{M_t} = 1$.

The set of all such uniformly balanced, uniformly bounded in $\|\cdot\|_{D_t}$ processes will be denoted $\D$.
\end{definition}

\begin{definition}
 A driver $g$ will be called \emph{uniformly balanced} if there exists a process $r\in \D$ such that for any $t, y, z,z'$ of appropriate dimension,
\[|g(t, y, z) - g(t, y, z')| \leq \|r_t(z-z')\|_{M_t}.\]
\end{definition}

\begin{lemma}
 A uniformly balanced driver is balanced.
\end{lemma}
\begin{proof}
 We can see that, for any $z,z'\in\mathbb{R}^\infty$,
\[\begin{split}\frac{|g(t, y, z)-g(t,y,z')|}{\|z-z'\|^2_{M_t}}
|(z-z') \Delta M_t|&\leq \frac{\|r_t(z-z')\|_{M_t}}{\|z-z'\|^2_{M_t}}
|(z-z') \Delta M_t|\\
& =\left\|r_t \frac{z-z'}{\|z-z'\|_{M_t}}\right\|_{M_t} \left|\frac{z-z'}{\|z-z'\|_{M_t}} \Delta M_t\right|
  \end{split}\]
Writing $u=\frac{z-z'}{\|z- z'\|_{M_t}}$,  the result is clear from Lemma \ref{lem:wellbalancedisbalanced}.
\end{proof}

\begin{definition}
Let $r\in\D$. We shall denote by $\E^r$ the nonlinear expectation given by $\E_g$ with $g(t,y,z)=\|r_tz\|_{M_t}$.

Similarly, we define $\E^{-r}$ to be the nonlinear expectation given by $\E_g$ with $g(t,y,z)=-\|r_tz\|_{M_t}$. 
\end{definition}

\begin{remark}
 As it is easy to show $\|r_t z\|^2_{M_t} \leq \sup_t(\|r_t\|_{D_t})^2  \|z\|^2_{M_t}$, the requirements for the existence of solutions to the BSDE are satisfied. As $r\in \D$, it is easy to show that $g(t,z)=\|r_tz\|_{M_t}$ is a uniformly balanced driver.
 
Note also that $\E^{r}$ is convex, positively homogenous and translation invariant, hence the up and downcrossing inequalities of Theorem \ref{thm:updowncrossing} apply.
\end{remark}

\begin{lemma}\label{lem:expgrowbound}
For any $Q$,
\[E[\E^r(Q|\F_t)^2]\leq E[Q^2] \exp\left((\sup_s\|r_s\|_{D_s}^2)(\mu_T-\mu_t)\right).\]
\end{lemma}
\begin{proof}
 Let $Y_t = \E^r(Q|\F_t)$. From the differentiation rule, we see that, for any predictable process $x_s>0$,
\[\begin{split}
   &E[Y_t^2] \\
&= E\left[Q^2 + 2\int_{]t,T]}\|r_sZ_s\|_{M_s} Y_{s-}d\mu_s - \int_{]t,T]} Z_s^2d\langle M\rangle_s - \sum_{t<s\leq T} \|r_s Z_s\|^2_{M_s}(\Delta\mu_s)^2\right]\\
& \leq E\left[Q^2 +\int_{]t,T]} x_s Y_{s-}^2 d\mu_s + \int_{]t,T]}(x_s^{-1}-\Delta\mu_s)\|r_sZ_s\|^2_{M_s} d\mu_s - \int_{]t,T]} Z_s^2d\langle M\rangle_s \right]
  \end{split}
\]
Setting $x_s^{-1} =\|r_s\|_{D_s}^{-2}+\Delta\mu_t$, we see from (\ref{eq:isometry}) that
\[\int_{]t,T]}(x_s^{-1}-\Delta\mu_s)\|r_sZ_s\|^2_{M_s} d\mu_s -\int_{]t,T]} Z_s^2d\langle M\rangle_s \leq 0\]
and $x_s\Delta\mu_s<1$. Hence we have
\[E[Y_t^2] \leq E[Q^2] +\int_{]t,T]} E[Y_{s-}^2] x_s d\mu_s\]
and an application of the Backward Gr\"onwall inequality (Lemma \ref{lem:groenwallinequality}) yields
\[E[\E^r(Q|\F_t)^2]\leq E[Q^2|\F_t] \mathfrak{E}(\tilde N;T)\mathfrak{E}(\tilde N;t)^{-1}\]
where $N_t = \int_{]0,t]}x_ud\mu_u$. Considering the continuous and discontinuous parts of $N$, we see that its right-jump-inverstion (Defintion \ref{defn:jumpinversion}) is $\tilde N_t = \int_{]0,t]}\|r_s\|_{D_s}^2 d\mu_s$, and hence
\[\begin{split}\mathfrak{E}(\tilde N;T)&=\mathfrak{E}(\tilde N;t)\exp\left(\int_{]t,T]}\|r_s\|_{D_s}^2 d\mu_t\right) \prod_{t<s\leq T} (1+\Delta\tilde N_s)e^{-\Delta\tilde N_s}\\
&\leq \mathfrak{E}(\tilde N;t)\exp((\mu_T-\mu_t)(\sup_s\|r_s\|_{D_s}^2)
  \end{split}
\]
yielding the result.
\end{proof}

\subsection{$\E^r$-dominated expectations}

\begin{definition}
For $r\in\D$, we say that a nonlinear expectation $\E$ is dominated by $\E^r$ if
\[\E(X+\eta) - \E(\eta) \leq \E^r(X)\]
for all $X, \eta \in L^2(\F_T)$.
\end{definition}

\begin{lemma}\label{lem:domrbound}
If $\E$ is dominated by $\E^r$, then
\[\E^{-r}(X) \leq \E(X+\eta) - \E(\eta) \leq \E^r(X)\]
for all $X, \eta \in L^2(\F_T)$.
\end{lemma}
\begin{proof} As noted in \cite{Coquet2002}, this is a simple consequence of the fact that $\E^{-r}(X) = -\E^r(-X)$.
\end{proof}

\begin{lemma}\label{lem:continuityofE}
If $\E$ is dominated by $\E^r$ for some $r\in \D$, then for all $\epsilon>0$, $\E$ is a continuous operator on $L^{2\vee(1+\epsilon)}(\F_T)$, in the sense that there exists $C_\epsilon$ such that
\[|\E(X) -\E(X')| \leq \|X-X'\|_{1+\epsilon}\]
\end{lemma}
\begin{proof}
This is a consequence of Lemmata \ref{lem:expectationbound} and \ref{lem:domrbound}.
\end{proof}

\begin{lemma}
For $\E$ an $\E^r$-dominated, translation invariant $\F$-expectation,
\[\E^{-r}(X|\F_t) \leq \E(X|\F_t)\leq \E^{r}(X|\F_t).\]
\end{lemma}

\begin{lemma}
Let $\E$ and $\E'$ be two translation invariant $\F$-expectations, both dominated by $\E^r$ for some $r\in \D$. If
\[\E(X)\leq \E'(X)\]
for all $X\in L^2(\F_T)$, then 
\[\E(X|\F_t)\leq \E'(X|\F_t)\]
up to evanescence.
\end{lemma}

\begin{theorem}\label{thm:cadlagmodn}
 Let $\E$ be an $\F$-expectation dominated by $\E^r$ for some $r\in\D$. Then an $\E$-martingale $Y \in S^2$ has a c\`adl\`ag modification.
\end{theorem}
\begin{proof}
 As $Y$ is an $\E$-martingale, we have that, for any $t\leq T$
\[Y_t = \E(Y_T|\F_t) \leq \E^r(Y_T|\F_t)\]
and so $Y$ is an $\E^r$-submartingale. As $\E^r$ is convex, translation invariant and positively homogeneous, we can apply Theorem \ref{thm:updowncrossing} to see that $Y$ almost surely admits left and right limits.

Define the c\`adl\`ag process $Y'_t := \lim_{s\downarrow t} Y_s = Y_{t+}$, this limit being almost surely well defined. As we assume the usual conditions, $Y'$ is adapted. For any $t\leq T$, any $A\in \F_t$, we have $Y'_t I_A = \lim_{s\downarrow t} Y_s I_A$, taking the limit in $L^2$ (which converges as $Y\in S^2$). From Lemma \ref{lem:continuityofE}, we see that $\E(Y'_tI_A) = \lim_{s\downarrow t} \E(Y_s I_A)$, but also, as $Y$ is an $\E$-martingale,
\[\E(Y_s I_A) = \E(\E(Y_s|\F_t)I_A) = \E(Y_t I_A)\]
and so $Y'_t = Y_t$ almost surely.
\end{proof}

\section{Doob-Meyer Decomposition for $g$-expectations}

We now show that, for a $g$-expectation $\E_g$, a Doob-Meyer decomposition holds. The method of proof is based on those in Peng \cite{Peng1999} (see also Royer \cite{Royer2006}). We begin with an $\E_g$-supermartingale $Y$ with $E[\sup_t (Y_t)^2]<\infty$. We wish to show that $Y$ can be written in the form
\[Y_t = Y_0 - \int_{]0,t]}g(u, Y_{u-}, Z_u) d\mu_u - A_t + \int_{]t,T]} Z_udM_u,\]
for some nondecreasing c\`adl\`ag process $A$ with $A_0=0$.

Similarly to \cite{Peng1999}, we shall use a sequence of penalised BSDEs. Consider the sequence of BSDEs with terminal values $Y_T^n=Y_T$, and drivers
\[f^n(t,y,z) = g(t,y, z) + n(Y_{t-}-y)^+.\]
The solutions of these BSDEs will be denoted\footnote{Note that this is a slight abuse of notation, as $Z^n$ here refers not to the $n$th component of $Z$, but the the $\mathbb{R}^\infty$ valued process which solves the BSDE with driver $f^n$. We shall not need to refer to individual components of $Z$ hereafter, and so this should not lead to confusion.}  $(Y^n, Z^n)$.

\begin{lemma}
The BSDEs with terminal values $Y_T$ and drivers $f^n$ have solutions $(Y^n, Z^n)$, which satisfy
\[\E(Y_T|\F_t) = Y^0_t \leq Y^n_t \leq Y^{n+1}_t \leq Y_t\]
and $Y^n_t\uparrow Y_t$ pointwise, up to evanescence. Furthermore $\{Y^n\}$ is a uniformly bounded set in $S^2$, and  $Y^n_{\cdot-} \to Y_{\cdot-}$ in $H^2$, that is,
\[E\left[\int_{]0,T]} \|Y^n_{t-}-Y_{t-}\|^2d\mu_t\right] \to 0.\]
\end{lemma}
\begin{proof}
As $g$ is firmly Lipschitz continuous, we have solutions for $f^0$ by Theorem \ref{thm:BSDEExist2}. For $n>0$, we can apply the same measure change argument as in \cite[Theorem 6.1]{Cohen2010} to assume without loss of generality that the Lipschitz constant of $g$ with respect to $y$ satisfies $c_t\Delta\mu_t<1-\epsilon$ for some $\epsilon>0$, and furthermore, $c>\epsilon^{-1}-1$. Hence we see that $f^n$ satisfies the requirements for Theorem \ref{thm:BSDEExist3}. Therefore these equations have solutions $(Y^n, Z^n)$.

By the comparison theorem (noting that $f^n$ is balanced as $g$ is balanced), we can see that $Y^n_t$ is nondecreasing in $n$ for all $t$, and that $Y^0_t=\E(Y_T|\F_t)$. Also if $Y^n_t >Y_t$, then by right continuity this must hold on some optional  interval $]\sigma, \tau]$, with $Y_\tau \geq Y^n_\tau$. However, on $]\sigma, \tau]$,  $Y^n_t = \E_g(Y^n_\tau|\F_t) \leq \E(Y_\tau|\F_t) \leq Y_t$ leading to a contradiction. Hence $Y^n_t \leq Y_t$ for all $n$, and all $t$. Therefore we have, for all $n$ and all $t$,
\[\E(Y_T|\F_t) = Y^0_t \leq Y^n_t \leq Y^{n+1}_t \leq Y_t.\]
Furthermore, suppose for some $\epsilon>0$, on some optional set $A$ nonempty with positive probability, we had $Y^{n}_t< Y_t-\epsilon$ for all $n$, all $t\in A$. Then $E[\int_{]0,T]} n (Y_{t-}-Y^n_{t-})^+ d\mu_t] \to \infty$, hence $Y^n_0 \to \infty$, which is a contradiction. Therefore, by continuity, $Y^n_t\uparrow Y_t$ except possibly on an evanescent set. By the dominated convergence theorem, it follows that $Y^n$ is a uniformly bounded set in $S^2$, and $Y^n_{\cdot-} \to Y_{\cdot-}$ in $H^2_\mu$.
\end{proof}

\begin{lemma}\label{lem:ZnAnbounded}
 Let $A^n_t=n\int_{]0,t]}(Y_{s-}-Y^n_{s-})^+d\mu_s$. Then there exists a constant $C$ independent of $n$ such that $E\left[\int_{]0,T]}(Z_t^n)^2d\langle M\rangle_t\right]<C$ and $E[(A^n_T)^2]<C$.
\end{lemma}
\begin{proof}
 From Ito's formula applied to $Y^n$, we see that,
\[\begin{split}&E[(Y^n_t)^2]+E\left[\int_{]t,T]} Z^n d\langle M\rangle_u\right] + E\left[\sum_{u\in ]t,T]}(g(u,Y^n_{u-},Z^n_u)\Delta\mu_u+\Delta A^n_u)^2\right] \\
&\qquad= E[Y_T^2] +2E\left[\int_{]t,T]}Y^n_{u-}(g(u,Y^n_{u-},Z^n_u)d\mu_u+dA^n_u)\right]
\end{split}\]
and hence,
\begin{equation}\label{eq:znbound1}
E\left[\int_{]t,T]} Z^n d\langle M\rangle_u\right]\leq E[Y_T^2] +2E\left[\int_{]t,T]}Y^n_{u-}g(u,Y^n_{u-},Z^n_u)d\mu_u\right]+2E\left[\int_{]t,T]}Y^n_{u-}(dA^n_u)\right]
\end{equation}
For $c$ the Lipschitz constant of $g$, we also have
\begin{equation}\label{eq:znbound2}
\begin{split}
&2E\left[\int_{]t,T]}Y^n_{u-}g(u,Y^n_{u-},Z^n_u)d\mu_u\right] \\
&\leq 4cE\left[\int_{]t,T]}(Y^n_{u-})^2d\mu_u\right] + (4c^{-1})E\left[\int_{]t,T]}(g(u,Y^n_{u-},Z^n_u))^2d\mu_u\right]\\
&\leq 4cE\left[\int_{]t,T]}(Y^n_{u-})^2d\mu_u\right] + (4c^{-1})E\left[\int_{]t,T]}(c(Y^n_{u-})^2 + c\|Z^n_{u-}\|^2_{M_u} + g(u,0,\mathbf{0})^2)d\mu_u\right]
\end{split}
\end{equation}
and
\begin{equation}\label{eq:znbound3}
\begin{split}
&2E\left[\int_{]t,T]}Y^n_{u-}(dA^n_u)\right]\leq 2E[A^n_T(\sup_{u}|Y^n_{u-}|)]\leq 2E[\sup_{u}(Y^n_{u-})^2]^{1/2}E[(A_T^n)^2]^{1/2} \\&\qquad\leq (16c\mu_T+8)E[\sup_{u}(Y^n_{u-})^2] + (16c\mu_T+8)^{-1}E[(A_T^n)^2]
\end{split}
\end{equation}
As $(Y^n)^2 \leq (Y^0)^2 + Y^2\in S^2$ and $E\left[\int_{]t,T]}\|Z^n_{u-}\|^2_{M_u} d\mu_u\right]  \leq E\left[\int_{]t,T]} Z^n d\langle M\rangle_u\right]$, combining (\ref{eq:znbound1}), (\ref{eq:znbound2}) and (\ref{eq:znbound3}), it follows that there is a constant $C_1$ independent of $n$ such that
\begin{equation}\label{eq:znbound4}
E\left[\int_{]t,T]} Z^n d\langle M\rangle_u\right]\leq C_1 + (8c\mu_T+4)^{-1}E[(A_T^n)^2].
\end{equation}

Furthermore, we also have
\[\begin{split}
A^n_T &= Y^n_0-Y^n_T -\int_{]0,T]}g(u,Y^n_{u-},Z^n_u)d\mu_u + \int_{]0,T]}Z_u^n dM_u\\
&\leq |Y_0|+|Y_T| +\int_{]0,T]}|g(u,Y^n_{u-},Z^n_u)|d\mu_u + \left|\int_{]0,T]}Z_u^n dM_u\right|
\end{split}\]
from which, expanding $(g(u,Y^n_{u-},Z^n_u))^2$ as in (\ref{eq:znbound2}), it follows that there exists a constant $C_2$ independent of $n$ such that
\begin{equation}\label{eq:znbound5}
\begin{split}E[(A^n_T)^2] &\leq 4E[(|Y_0|+|Y_T|)^2] + 4\mu_TE\left[\int_{]0,T]}(g(u,Y^n_{u-},Z^n_u))^2d\mu_u\right] +2E\left[\left(\int_{]0,T]}Z_u^n dM_u\right)^2\right]\\
&\leq C_2 + (4c\mu_T+2)E\left[\int_{]t,T]} Z^n d\langle M\rangle_u\right].\\
\end{split}\end{equation}
Combining (\ref{eq:znbound4}) and (\ref{eq:znbound5}) yields the result. 
\end{proof}

We can now prove the convergence of our solutions. Unlike in \cite{Peng1999} and \cite{Royer2006}, due to the use of left-limits in the BSDE, we are able to prove the strong convergence of $Z^n$ in $L^2$, rather than only in $L^p$ for $p<2$.
\begin{theorem}\label{thm:DMdecompgen}
A c\`adl\`ag $\E_g$-supermartingale $Y$ has a representation of the form
\[Y_t = Y_0 - \int_{]0,t]}g(u, Y_{u-}, Z_u) d\mu_u - A_t + \int_{]0,t]} Z_udM_u,\]
where $Z$ is the strong limit of $Z^n$ in $H^2_M$ and $A$ is a c\`adl\`ag increasing process.
\end{theorem}
\begin{proof}
By Lemma \ref{lem:ZnAnbounded}, we know that $\{Z^n\}_{n\in\mathbb{N}}$ is weakly compact in $H^2_{M}$, and, defining $g^n_t:=g(t, Y^n_{t-}, Z^n_t)$, we see $\{g^n\}$ is bounded and hence weakly compact in $H^2_{\mu}$. Therefore, by extracting subsequences, we have the existence of weak limits $Z^n\rightharpoonup Z$ and $g^n\rightharpoonup g^\infty$. For any stopping time $\tau\leq T$, we also then have the weak convergence of the integrals $\int_{]0,\tau]}Z^n_udM_u$ and $\int_{]0,\tau]}g^n_ud\mu_u$ in $L^2(\F_T)$.  As
\[A^n_t = Y_0^n-Y_t^n - \int_{]0,t]} g^n_u d\mu_u + \int_{]0,t]} Z^n_udMu\]
we also have the existence of a weak $L^2$-limit 
\[A^n_t \rightharpoonup A_t = Y_0-Y_t - \int_{]0,t]} g^\infty_u d\mu_u + \int_{]0,t]} Z_udM_u\]
and clearly, $A$ is a nondecreasing process with $A_T\in L^2(\F_T)$. By a result of Peng \cite[Lemma 2.2]{Peng1999}, $A$ is c\`adl\`ag. As $Y$ is given, we see that $Z$ is uniquely defined, and hence the sequence $\{Z^n\}$ (rather than a subsequence) must weakly converge.

We now write $\delta_n Y = Y-Y^n$, $\delta_n Z = Z-Z^n$, $\delta_n g = g^\infty-g^n$ and $\delta_n A = A-A^n$. Considering the dynamics of $(\delta_n Y)^2$, from It\^o's formula we have
\[\begin{split}0&=E[\delta_n Y^2_T]\\
&=E[\delta_n Y^2_0] -2E\left[\int_{]0,T]} (\delta_n Y)_{u-} ((\delta_n g_u) d\mu + d(\delta_n A)_u)\right]\\
&\qquad + E\left[\int_{]0,T]} (\delta_n Z_u)^2d\langle M\rangle_u\right] + E\left[\sum_{u\in]0,T]} ((\delta_n g_u )\Delta \mu_u + \Delta (\delta_nA)_u)^2\right]\end{split}\]
from which we obtain
\[E\left[\int_{]0,T]} (\delta_n Z_u)^2d\langle M\rangle_u\right] \leq 2E\left[\int_{]0,T]} (\delta_n Y)_{u-} ((\delta_n g_u) d\mu + d(\delta_n A)_u)\right].\]
We then see that, by the Cauchy-Schwartz inequality, for $C$ a bound on $\delta_n g$ in $H^2_\mu$,
\[\begin{split}
E\left[\int_{]0,T]} (\delta_n Y)_{u-}(\delta_n g_u) d\mu \right]
&\leq E\left[\int_{]0,T]} (\delta_n Y)_{u-}^2 \right]^{1/2} E\left[\int_{]0,T]} (\delta_n g_u)^2 d\mu \right]^{1/2}\\
&\leq C\cdot E\left[\int_{]0,T]} (\delta_n Y)_{u-}^2 \right]^{1/2} \\
&\rightarrow 0.
\end{split}\]
Also
\[\begin{split}E\left[\int_{]0,T]} (\delta_n Y)_{u-} d(\delta_n A)_u\right]&=E\left[\int_{]0,T]} (\delta_n Y)_{u-} dA_u\right] -E\left[\int_{]0,T]} (\delta_n Y)_{u-} dA^n_u\right]\\
&\leq E\left[\int_{]0,T]} (\delta_n Y)_{u-} dA_u\right]\\
&\leq E[A_T\sup_u (\delta_0 Y)_u]\\
&\leq E[A_T^2]+E[\sup_u (\delta_0 Y)^2_u] <\infty
\end{split}\]
and so, by the Dominated convergence theorem, 
\[E\left[\int_{]0,T]} (\delta_n Y)_{u-} d(\delta_n A)_u\right]\leq E\left[\int_{]0,T]} (\delta_n Y)_{u-} dA_u\right] \to 0.\] 
Hence we see that,
\[E\left[\int_{]0,T]} (\delta_n Z_u)^2d\langle M\rangle_u\right] \to 0.\]
Given this strong convergence, it is clear that $g^n\to g^{\infty}$ strongly in $H^2_\mu$, and that $g^\infty_t = g(t, Y_{t-}, Z_t)$ $\mathbb{P}\times\mu$-a.e., yielding the desired representation. 
\end{proof}

To compare this with the classical Doob-Meyer decomposition, we have the following corollary.
\begin{corollary}
Consider $\E_g$ a $g$-expectation, where $g(u, z)$ does not depend on $y$ (and hence, $\E_g$ is translation invariant). Then a c\`adl\`ag $\E_g$-supermartingale $Y$ in $S^2$ has a decomposition $Y=Y_0+M-A$, where $A$ is a nondecreasing adapted c\`adl\`ag process with $A_T\in L^2(\F_T)$, and $M$ is a c\`adl\`ag $\E_g$-martingale in $S^2$ with $M_0=0$.
\end{corollary}
\begin{proof}
From Theorem \ref{thm:DMdecompgen}, we have the representation
\[Y_t = Y_0 - \int_{]0,t]}g(u, Z_u) d\mu_u - A_t + \int_{]0,t]} Z_udM_u,\]
and note that
\[\begin{split}
M_t &= - \int_{]0,t]}g(u, Z_u) d\mu_u+ \int_{]t,T]} Z_udM_u \\
&= M_T + \int_{]t,T]}g(u, Z_u) d\mu_u- \int_{]t,T]} Z_udM_u =\E_g(M_T|\F_t)
\end{split}\]
is a $g$-martingale.
\end{proof}

We can now show that $\E^r$-domination implies that the drift must be $\mu$-absolutely continuous.

\begin{theorem}\label{thm:driftabscont}
 Let $\E$ be an $\F$-expectation, $\E^r$-dominated for some $r\in D^\infty$. Let $Y$ be a c\`adl\`ag $\E$-martingale. Then there exist unique predictable processes $g\in H^2_\mu$, $Z\in H^2_M$ such that
\[Y_T = Y_t - \int_{]t,T]} g_ud\mu_u + \int_{]t,T]}Z_udM_u\]
up to indistinguishability. These processes satisfy $|g_u| \leq \|r_uZ_u\|_{M_t}$.
\end{theorem}

\begin{proof}
 As $\E$ is $\E^r$-dominated, we know that
\[\E^{-r}(Y_T|\F_t) \leq Y_t \leq \E^r(Y_T|\F_t)=-\E^{-r}(-Y_T|\F_t) ,\]
and so both $Y$ and $-Y$ are $\E^{-r}$-supermartingales. From the nonlinear Doob-Meyer decomposition (Theorem \ref{thm:DMdecompgen}), we can find nondecreasing c\`adl\`ag processes $A^r$, $A^{-r}$ and processes $Z^r, Z^{-r} \in H^2_M$ such that
\begin{equation}\label{eq:dmEmartpair}
 \begin{split}
  Y_t &= Y_0 + \int_{]0,t]} \|r_uZ^{-r}_u\|_{M_u} d\mu + \int_{]0,t]} Z^{-r}_udM_u -A^{-r}_t\\
  -Y_t &= -Y_0 + \int_{]0,t]} \|r_uZ^{r}_u\|_{M_u} d\mu + \int_{]0,t]} Z^r_udM_u -A^r_t.
  \end{split}
\end{equation}
As $Y$ is a special semimartingale, its canonical decomposition (into martingale and predictable finite-variation components) is unique (see \cite[Def 4.22]{Jacod2003}). Hence we have $\int_{]0,t]} Z^{-r}dM = -\int_{]0,t]} Z^rdM$ up to indistinguishability, and furthermore $Z^{-r} = -Z^r$ in $H^2_M$.
Taking the sum of the two equations in (\ref{eq:dmEmartpair}), we then have
\[0 = 2\int_{]0,t]} \|r_uZ^{r}_u\|_{M_u} d\mu-A^{-r}_t-A^r_t.\]
Differentiating yields
\[d(A^r+ A^{-r})_u = 2\|r_uZ^{r}_u\|_{M_u} d\mu\]
and, as both $A^r$ and $A^{-r}$ are nondecreasing, we see that they are both absolutely continuous with respect to $\mu$. Therefore, as $A^{-r}_T\in L^2(\F_T)$,  we can write $dA^{-r}_t = a^{-r}_t d\mu$ for some $a^{-r} \in H^2_\mu$.  Defining $g_u := -\|rZ^{-r}_u\|_{M_u} + a^{-r}_u$, we have
\[Y_t = Y_0 - \int_{]0,t]} g_u d\mu + \int_{]0,t]} Z^{-r}_udM_u.\]
This $g$ is unique among predictable processes in $H^2_\mu$, again by the uniqueness of the canonical decomposition of a special semimartingale. Furthermore, as $A^{-r}$ and $A^r$ are nondecreasing, we have that $0\leq a^{-r}\leq 2\|r_uZ^{r}_u\|_{M_u}$, and so $|g_u| \leq \|r_uZ_u\|_{M_t}$.
\end{proof}

\begin{theorem}\label{thm:induceddriftisunifbalanced}
 Let $\E$ be as in Theorem \ref{thm:driftabscont}, and $Y$ and $Y'$ be two c\`adl\`ag $\E$-martingales, with associated processes $g, g'$ and $Z, Z'$. Then
\[|g_t-g'_t| \leq \|r_t(Z_t-Z_t')\|_{M_t}\]
up to evansescence.
\end{theorem}
\begin{proof}
 As all of $Y, -Y, Y'$ and $-Y'$ are $\E^{-r}$-supermartingales, by Lemma \ref{lem:sumsupissup} we know that $\delta Y:=Y-Y'$ and $-\delta Y$ are both $\E^{-r}$-supermartingales. By precisely the same argument as in Theorem \ref{thm:driftabscont}, we can find predictable processes $g^\delta\in H^2_\mu$, $Z^\delta \in H^2_M$ such that
\[\delta Y_t = \delta Y_0 - \int_{]0,t]} g^\delta d\mu + \int_{]0,t]} Z^\delta_u dM\]
and $|g^\delta_t| \leq \|r_tZ^\delta_t\|_{M_t}$ up to evanescence. However, we also have
\[\delta Y_t = \delta Y_0 - \int_{]0,t]} (g_u-g'_u)d\mu + \int_{]0,t]}(Z_u-Z'_u)dM\]
and uniqueness of the canonical decomposition of $\delta Y_t$ yields
\[|g_t-g'_t|=|g^\delta_t| \leq \|r_tZ^\delta_t\|_{M_t} = \|r_t(Z_u-Z'_u)\|_{M_t}.\]
\end{proof}

\section{$\E^r$-dominated Doob-Meyer decomposition}
We shall need to extend our decomposition to the case where $\E$ is $\E^r$-dominated for some $r\in\D$, but where we do not know \emph{a priori} that it is a $g$-expectation.

We need the following generalisation of our existence result. A more general result than this is possible (where $n(Y_{t-}-y_{t-})^+$ is replaced by an appropriately Lipschitz function with sufficiently bounded upward jumps). This is, however, largely pointless given the representation we shall prove further on (Theorem \ref{thm:representation}), which implies these results are equivalently given by Theorem \ref{thm:DMdecompgen}.

\begin{theorem}\label{thm:Erdomdoobmeyer}
 Consider $\E$ any translation invariant $\F$-expectation, $\E^r$-dominated for some $r\in\D$. For  any $Q\in L^2(\F_T)$, any c\`adl\`ag $\E$-supermartingale $Y$ in $H^2_\mu$ with $Y_T=Q$, the equation
\[Y^n_t= \E\left.\left(Q+n\int_{]t,T]}(Y_{u-} - Y^n_{u-})^+d\mu\right|\F_t\right)\]
has a unique c\`adl\`ag solution in $H^2_\mu$.
\end{theorem}

\begin{proof}
 Our approach is similar to that in Theorem \ref{thm:BSDEExist3}. For any $s<t$, any $Q'\in L^2(\F_{t-})$, define a mapping
\[\Phi_{]s,t[}^{Q'}: H^2_\mu \to H^2_\mu,\qquad y \mapsto \E\left.\left(Q'+n\int_{]s,t[}(Y_{u-} - y_{u-})^+d\mu\right|\F_t\right).\]
For any two approximations $y,y'\in H^2_\mu$, define $\delta y = y-y'$ and $\delta\Phi(y) = \Phi_{]s,t[}^{Q'}(y)-\Phi_{]s,t[}^{Q'}(y')$. Then as $\E$ is $\E^r$-dominated, and $r$ is assumed to be bounded (as it is uniformly balanced), it is easy to show (see \cite[Lemma 6.1]{Coquet2002} and use Lemma \ref{lem:expgrowbound})
\[\begin{split}
   E[(\delta \Phi(y))^2]
&\leq E\left[\left(n \E^r\left.\left(\int_{]s,t[}|\delta y|d\mu\right|\F_t\right)\right)^2\right]\\
&\leq n^2 e^{(\sup_t\|r_t\|_{D_t}^2)(\mu_{t-}-\mu_s)} E\left[\left(\int_{]s,t[}|\delta y|d\mu\right)^2\right]\\
&\leq n^2 e^{\sup_t\|r_t\|_{D_t}^2\mu_T}(\mu_{t-}-\mu_s)E\left[\int_{]s,t[}|\delta y|^2d\mu\right].
\end{split}\]
As $\mu$ is summable, using the result of \cite[Lemma 6.1]{Cohen2010}, we can find a finite set $\{0=t_1, t_2,..., t_m=T\}$ where $n^2 e^{\sup_t\|r_t\|_{D_t}^2\mu_T}(\mu_{t_{i+1}-}-\mu_{t_i})<1$ for all $i$. Hence we have a contraction on each of the subintervals $]t_i, t_{i+1}[$. Therefore, for any $Y^n_{t_{i+1}-} = Q'\in L^2(\F_{t_{i+1}-})$, we can solve our equation uniquely back to time $t_i$.

At each $t_i$, we shall solve the equation directly. Suppose we have a solution $Y^n_u$ for all $u\geq t_i$. In particular, we have the value $Y^n_{t_i} \in L^2(\F_{t_i})$. Then we have the equation
\[Y^n_{t_i-} = \E\left.\left(Y_{t_i}+n(Y_{t_i-} - Y^n_{t_i-})^+\Delta\mu\right|\F_{t-}\right)\]
which, by translation invariance of $\E$, gives
\[Y^n_{t_i-} = \left(\frac{1}{1+n\Delta \mu_{t_i}} \E\left.\left(Y^n_{t_i} + n\Delta\mu_t Y^n_{t_i-}\right|\F_{t-}\right)\right)\wedge \E(Y^n_{t_i}|\F_{t-})\]
Note as $n\Delta\mu_t>0$, $Y^n_{t_i-}$ is clearly in $L^2(\F_{t-})$. Therefore, at each time $t_i$, we can take any $Y_{t_i}\in L^2(\F_{t_{i+1}})$, and obtain a unique value $Y_{t_i-}\in L^2(\F_{t_{i+1}-})$.

 Using backward induction and alternating between the contraction mapping approach and the direct approach yields a unique solution. It is then straightforward to verify (as in Theorem \ref{thm:BSDEExist3}) that this solution is c\`adl\`ag and in $H^2_\mu$.
\end{proof}

\begin{lemma}
 For $Y, Y^n$ as in Theorem \ref{thm:Erdomdoobmeyer},
\[\E(Q|\F_t) = Y^0 \leq Y^n \leq Y^{n+1} \leq Y.\]
\end{lemma}
\begin{proof}
That $Y^n\geq Y^0 = \E(Q|\F_t)$ is easy from the monotonicity of $\E$.

Suppose $Y^n_t\geq Y^{n+1}_t$ with positive probability. By right continuity, there exists an optional interval $A=]\sigma, \tau]$, nonempty with positive probability, such that $Y^n_t\geq Y^{n+1}_t$ on $]\sigma, \tau[$ and $Y^n_\tau\geq Y^{n+1}_\tau$. On $A$, note that $(Y-Y^n)^+ \leq (Y-Y^{n+1})^+$, and hence for any $t\in A$,
\[\begin{split}
   I_{A} Y_t^n
&= \E\left.\left(I_{A}Y^n_{\tau} + \int_{]t,\tau]} n I_{A}(Y_{u-}-Y^n_{u-})^+d\mu_u\right|\F_t\right)\\
&\leq \E\left.\left(I_{A}Y^{n+1}_{\tau} + \int_{]t,\tau]} (n+1) I_{A}(Y_{u-}-Y^{n+1}_{u-})^+d\mu_u\right|\F_t\right)\\
&=   I_{A} Y_t^{n+1}
  \end{split}\]
which gives a contradiction. Hence $Y^n\leq Y^{n+1}$. A similar argument applies with $Y^{n+1}$ replaced by $Y$.
\end{proof}

\begin{lemma}\label{lem:decomErdomsuper}
 For $Y^n$ as in Theorem \ref{thm:Erdomdoobmeyer}, $Y^n$ has a representation
\[Y_t^n = Y^n_0 -\int_{]0,t]} g^n d\mu - A^n_t + \int_{]0,t]}Z^ndM\]
for some $g^n\in H^2_\mu$, $Z^n\in H^2_M$ and  $A^n_t$ nondecreasing, predictable and c\`adl\`ag with $A_0=0$ and $A_T\in L^2(\F_T)$. Furthermore, $|g_u^n| \leq \|r_u Z^n_u\|_{M_u}$, and there exists a constant $C$ independent of $n$ such that $E[(A^n_T)^2]<C$ and $E[\int_{]0,T]} (Z^n)^2_u d\langle M\rangle_u]<C$.
\end{lemma}
\begin{proof}
Define $A^n_t = \int_{]0,t]} n (Y_u- Y^n_u)^+d\mu$. As $Y^n + \int_{]0,t]} n (Y_u- Y^n_u)^+d\mu$ is a $\E$-martingale, we have from Theorem \ref{thm:driftabscont} the existence of $g^n$ and $Z^n$ with the required inequality between them.

For the required bound on $E[(A^n_T)^2]$ and $E[\int_{]0,T]} (Z^n)^2_u d\langle M\rangle_u]$, as $|g^n_t|<\|r_t Z^n_t\|_{M_t}$, where $r\in D^\infty$, we can precisely repeat the argument of Lemma \ref{lem:ZnAnbounded}.
\end{proof}

\begin{theorem}\label{thm:DMdecompErdom}
 Let $\E$ be a translation invariant $\F$-expectation, which is $\E^r$-dominated for some $r$. A c\`adl\`ag $\E$-supermartingale $Y$ has a representation of the form
\[Y_t +A_t = \E(Y_T + A_T|\F_t)\]
where $A$ is a nondecreasing, predictable and c\`adl\`ag process with $A_T\in L^2(\F_T)$.
\end{theorem}

\begin{proof}
 As in the proof of Theorem \ref{thm:DMdecompgen}, we see that the $A^n$ and $Z^n$ terms constructed in Lemma \ref{lem:decomErdomsuper} are uniformly bounded, and so must weakly converge. As $|g^n_u|\leq \|r_uZ^n_u\|_{M_u}$, we can again see that the argument of Theorem \ref{thm:DMdecompgen} will hold, and so $Z^n$ converges strongly in $H^2_M$. Therefore $g^n$ converges strongly in $H^2_\mu$, and hence $A^n_t$ converges strongly in $L^2(\F_t)$. By Lemma \ref{lem:continuityofE}, we can pass to the $L^2$-limit in the equation $Y_t^n +A_t^n = \E(Y_T^n + A_T^n|\F_t)$, and the theorem is proven.
\end{proof}

\section{Representation as a $g$-expectation}

We can now prove our main result, that any translation invariant $\F$-expectation which is $\E^r$-dominated for some $r\in\D$, must be a $g$-expectation.

\begin{theorem}\label{thm:representation}
 Consider a translation invariant $\F$-expectation $\E$, which is $\E^r$-dominated for some $r\in \D$. Then there exists a unique function $g:\Omega \times [0,T] \times \mathbb{R}^\infty \to \mathbb{R}$ satisfying $E[\int_{]0,T]}(g(t,\mathbf{0}))^2d\mu]<\infty$ and $g$ is uniformly balanced (and hence Lipschitz), such that
\[\E(Q) = \E_g(Q)\]
for all $Q\in L^2(\F_T)$. Furthermore, $g(t, \mathbf{0}) =0$ for $\mu$-almost all $t$.
\end{theorem}
\begin{proof}
 For each $z\in \mathbb{R}^\infty$, we consider the forward equation
\[dY^z = -\|r_t z\|_{M_t} d\mu_t + zdM;\qquad Y^z_0=0.\]
 We then see that $Y^z$ is an $\E^r$-martingale, and hence an $\E$-supermartingale.

From Theorem \ref{thm:DMdecompErdom}, there exists a nondecreasing, predictable and c\`adl\`ag process $A^z$ with $A^z_0=0$ and $A^z_T\in L^2(\F_T)$ such that
\[Y^z_t+ A^z_t = \E(Y^z_T + A^z_T|\F_t).\]
By Theorem \ref{thm:driftabscont}, there is a unique $g(z; \cdot):\Omega\times[0,T]\to\mathbb{R}$ predictable such that
\[Y^z_t + A^z_t = Y^z_T + A^z_T + \int_{]t,T]} g(z;u)d\mu - \int_{]t,T]} Z^z_u dM_u.\]
and $|g(z;t)| \leq \|r_t Z^z\|_{M_t}$.

As we also know
\[Y^z_t = Y^z_T + \int_{]t,T]} \|r_t z\|_{M_t}d\mu - \int_{]t,T]}z dM_u\]
we see that
\[A^z_t \equiv \|r_t z\|_{M_t}- \int_{]0,t]} g(z;u)d\mu,\quad Z^z \equiv z.\]
In particular, this implies $|g(z;t)| \leq \|r_t z\|_{M_t}$. From Theorem \ref{thm:induceddriftisunifbalanced}, we also see that for any $z,z'\in\mathbb{R}^\infty$,  $|g(z;t)- g(z';t)| \leq \|r_t (z-z')\|_{M_t}$. Hence, for each $t$, $g(\cdot;t)$ is uniformly Lipschitz continuous and uniformly balanced, as a function of $z$.

We can see that, for any $0\leq r\leq t \leq T$,
\[Y^z_t + A^z_t = Y^z_r + A^z_r - \int_{]r,t]} g(z;u)d\mu + \int_{]r,t]} z dM_u.\]
Because of translation invariance, we have
\[\E\left.\left(- \int_{]r,t]} g(z;u)d\mu + \int_{]r,t]} z dM_u\right|\F_r\right)=0.\]
Let $\{A_i\}_{i=1}^N \subset \F_r$  be a partition of $\Omega$, and let $z_i\in \mathbb{R}^\infty$. From Lemma \ref{lem:regularityofEt}, and the fact $g(0,t)\equiv 0$, it follows that
\[\E\left.\left(- \int_{]r,t]} g\left(\sum_i I_{A_i}z_i;u\right)d\mu + \int_{]r,t]} \left(\sum_i I_{A_i}z_i\right) dM_u\right|\F_r\right)=0.\]
Hence, by the continuity of $\E$ given in Lemma \ref{lem:continuityofE} and the fact that $g$ is Lipschitz in $z$, we have, for any $Z\in H^2_M$,
\[\E\left.\left(- \int_{]r,t]} g\left(Z_u;u\right)d\mu + \int_{]r,t]} Z_u dM_u\right|\F_r\right)=0.\]

For any $Q\in L^2(\F_T)$, now solve the BSDE with driver $g$. As $g$ is Lipschitz, this has a unique solution $(Y, Z)$, and by the definition of $g$-expectation, $\E_g(Q)=Y_0$. On the other hand, we also have
\[\begin{split}\E(Q) &= \E\left(Y_0- \int_{]0,T]} g\left(Z_u;u\right)d\mu + \int_{]0,T]} Z_u dM_u\right) \\
   &= Y_0 + \E\left(- \int_{]0,T]} g\left(Z_u;u\right)d\mu + \int_{]0,T]} Z_u dM_u\right)=Y_0
  \end{split}
\]
and so $\E_g(Q)=Y_0 = \E(Q)$ for all $Q\in L^2(\F_T)$.
\end{proof}

\section{Conclusion}
We have extended the results of \cite{Coquet2002} and \cite{Peng1999} to a general setting. This directly answers the question raised by Remark 7.1 of \cite{Coquet2002}; we have given a nonlinear Doob-Meyer decomposition theorem for $g$-expectations, and have shown that every $\F$-expectation satisfying a dominance relation can be expressed as a $g$-expectation. Our only assumption on the probability space is that $L^2(\F_T)$ is separable.

The exact nature of this dominance relation is quite interesting in this context. One can think of the dominance relation in \cite{Coquet2002} as being needed to guarantee that the induced driver of the BSDE exists, and is Lipschitz continuous. Our assumption guarantees both these properties, and furthermore that the driver can be integrated with respect to the (arbitrary) Stieltjes measure $\mu$, and that it satisfies the conditions to be uniformly balanced, and so a comparison theorem will hold. Neither of these properties appears in \cite{Coquet2002}, as in they assume that $\mu$ is always Lebesgue measure (a reasonable assumption, as all martingales have absolutely continuous quadratic variation), and all martingales are continuous (so the comparison theorem holds automatically). However, if our filtration is generated by finitely many Brownian motions, as in \cite{Coquet2002}, then our result corresponds precisely to theirs. Furthermore, our result will also encompass the case of a filtration generated by countably many independent Brownian motions.

As $\D$ contains a wide range of processes, our assumption that $\E$ is $\E^r$-dominated for some $r\in \D$ has particular implications for those cases where the BSDE can be written in the form  (c.f. \cite{El1997a})
\[dY_t = -g(t,  Z_t)d\mu_t + Z_tdM'_t + dN_t\]
for some finite-dimensional martingale $M'$, where $N$ is a martingale orthogonal to $M'$. From the perspective of the Davis-Varaiya martingale representation theorem, this means that the BSDE driver looks only at a finite dimensional subspace of the space of $S^2$-martingales. Looking from the perspective of the $\F$-expectation, this is equivalent to stating that $\E(Q+N)=\E(Q)$ for any $Q$ and any martingale $N$ orthogonal to $M'$ and with $N_0=0$. In this context, if $\E$ is $\E^r$-dominated for some $r\in\D$, we can find a \emph{degenerate} matrix $r'\in \D$ such that $\E$ is $\E^{r'}$-dominated, and the representation will follow.

If we compare our results with the L\'evy case considered by Royer \cite{Royer2006}, we see that our condition `$g$ is uniformly balanced' is equivalent to her `assumption $A_\gamma$'. Royer shows that assumption $A_\gamma$ is satisfied by the BSDEs generated by nonlinear expectations, and we similarly show that the induced $g$ is uniformly balanced.

If we compare with earlier results in discrete time (\cite{Cohen2008c}, \cite{Cohen2009a}), we see that we have again shown an equivalence between BSDE solutions and translation invariant nonlinear expectations. Unlike in discrete time, we require the further assumption of $\E^r$-domination to ensure that the continuous-time generator is adequately Lipschitz continuous, and so our results lack the complete generality of those in discrete time.

Further work on this area may allow us to extend away from the assumption of translation invariance (see \cite{Cohen2009a} in discrete time), and towards quadratic BSDEs (see \cite{Hu2008} in the Brownian case). A further extension would also be to allow $\mu$ to be a stochastic finite-variation process.  These results will require further extension of the existence results of BSDEs in general filtrations.

\bibliographystyle{plain}
\bibliography{../RiskPapers/General}
\end{document}